\documentclass[12pt,reqno,a4wide]{amsart}

\usepackage{tikz, amssymb} 
\usepackage{url}
\usepackage{hyperref, adjustbox, float}
\usepackage{pstricks,pst-node,pst-tree}
\usetikzlibrary{decorations.pathreplacing}
\usetikzlibrary{fit}
\newcommand\addvmargin[1]{
  \node[fit=(current bounding box),inner ysep=#1,inner xsep=0]{};
}

\allowdisplaybreaks

\newtheorem{definition}{Definition}

\begin{document}

	\bibliographystyle{plain}

	\title[A bijection between ternary trees and a subclass of Motzkin paths]{A bijection between ternary trees and a subclass of Motzkin paths}

	\author[Helmut Prodinger]{Helmut Prodinger}
	\address{Department of Mathematics, University of Stellenbosch 7602,
		Stellenbosch, South Africa}
	\email{hproding@sun.ac.za}

	\author[Sarah J. Selkirk]{Sarah J. Selkirk}
	\address{Department of Mathematics, University of Stellenbosch 7602,
		Stellenbosch, South Africa}
	\email{sjselkirk@sun.ac.za}
	\thanks{The financial assistance of the National Research Foundation (NRF) towards this research is hereby acknowledged. Opinions expressed and conclusions arrived at, are those of the authors and are not necessarily to be attributed to the NRF}

	
	\begin{abstract}
		A bijection between ternary trees with $n$ nodes and a subclass of Motzkin paths of length $3n$ is given. This bijection can then be generalized to $t$-ary trees. 
	\end{abstract}
	

	\maketitle

	\section{Introduction}

A recent question in the International Mathematics Competition proposed by Petrov and Vershik \cite{IMC} counts the number of allowed paths from $(0, 0, 0)$ to $(n, n, n)$ of a frog that makes steps of length one along the lattice 
\begin{equation*}
\Omega = \{(x, y, z)\in \mathbb{Z}^{3} \mid 0 \leq z \leq y \leq x \leq y+1\}
\end{equation*}
in exactly $3n$ moves. 

Clearly there are $n$ steps in each of the three possible directions, and we model each step as follows:
\begin{center}
\begin{tabular}{| c | c | c | c |}
\hline
$(x, y, z)$  & $(0, 0, 1)$ & $(0, 1, 0)$ & $(1, 0, 0)$ \\ 
\hline
Step  & \begin{tikzpicture}[scale = 0.3, line width = 0.3mm]
\coordinate (aux) at (0,0);
\foreach \i in {-1}
    \draw[line cap = round] (aux)--++(1,\i) coordinate (aux);
\addvmargin{1mm}
\end{tikzpicture}
& \begin{tikzpicture}[scale = 0.3, line width = 0.3mm]
\coordinate (aux) at (0,0);
\foreach \i in {1}
    \draw[line cap = round] (aux)--++(1,\i) coordinate (aux);
\addvmargin{1mm}
\end{tikzpicture}
& \begin{tikzpicture}[scale = 0.3, line width = 0.3mm]
\coordinate (aux) at (0,0);
\foreach \i in {0}
    \draw[line cap = round] (aux)--++(1,\i) coordinate (aux);
\addvmargin{1mm}
\end{tikzpicture} \\
\hline
\end{tabular}
\end{center}

This along with the restriction $0 \leq z \leq y \leq x \leq y+1$, gives rise to the subclass of Motzkin paths defined below. 
\begin{definition}
An \textbf{S-Motzkin path} is a Motzkin path with $n$ of each type of step such that the following conditions hold
\begin{itemize}
\item The initial step must be 
\begin{tikzpicture}[scale = 0.3, line width = 0.3mm]
\coordinate (aux) at (0,0);
\foreach \i in {0}
    \draw[line cap = round] (aux)--++(1,\i) coordinate (aux);
\end{tikzpicture}\,, 
\item between every two 
\begin{tikzpicture}[scale = 0.3, line width = 0.3mm]
\coordinate (aux) at (0,0);
\foreach \i in {0}
    \draw[line cap = round] (aux)--++(1,\i) coordinate (aux);
\end{tikzpicture} 
there is exactly one 
\begin{tikzpicture}[scale = 0.3, line width = 0.3mm]
\coordinate (aux) at (0,0);
\foreach \i in {1}
    \draw[line cap = round] (aux)--++(1,\i) coordinate (aux);
\end{tikzpicture}\,, 
\item the $k$-th occurring 
\begin{tikzpicture}[scale = 0.3, line width = 0.3mm]
\coordinate (aux) at (0,0);
\foreach \i in {-1}
    \draw[line cap = round] (aux)--++(1,\i) coordinate (aux);
\end{tikzpicture}
must occur after at least $k$ pairs of 
\begin{tikzpicture}[scale = 0.3, line width = 0.3mm]
\coordinate (aux) at (0,0);
\foreach \i in {0}
    \draw[line cap = round] (aux)--++(1,\i) coordinate (aux);
\end{tikzpicture}
and 
\begin{tikzpicture}[scale = 0.3, line width = 0.3mm]
\coordinate (aux) at (0,0);
\foreach \i in {1}
    \draw[line cap = round] (aux)--++(1,\i) coordinate (aux);
\end{tikzpicture}\,.
\end{itemize}
\end{definition}

The total number of such paths is $\frac{1}{2n+1}\binom{3n}{n}$ which is equal to the number of ternary trees with $n$ nodes \cite{Knuth}. We first provide a mapping from S-Motzkin paths to ternary trees, and then provide the inverse mapping, thus showing that S-Motzkin paths are bijective to ternary trees as well as other combinatorial objects found in \cite{Gu, Panholz, Prod}. For completeness, an instructive example is given along with a table for $n=3$. 

\section{Bijection}

\subsection{S-Motzkin paths to ternary trees}

We define $\varnothing$ to be the empty path. For an arbitrary S-Motzkin path $\mathcal{M}$, the canonical decomposition is
\begin{equation*}
\Phi(\mathcal{M}) = \left( \mathcal{A},\, \mathcal{B},\, \mathcal{C}\right),
\end{equation*}
where $\mathcal{A}$, $\mathcal{B}$, and $\mathcal{C}$ represent paths at the left, middle, and right subtrees respectively. Furthermore,  
\begin{itemize}
\item $\mathcal{C}$ is the path from the penultimate return of the path to the last return, with the initial and last steps removed,
\item $\mathcal{A}$ is the path from $y$ to $x$ (not including $x$), where $x$ is the first \begin{tikzpicture}[scale = 0.3, line width = 0.3mm]
\coordinate (aux) at (0,0);
\foreach \i in {0}
    \draw[line cap = round] (aux)--++(1,\i) coordinate (aux);
\end{tikzpicture} 
to the left of $\mathcal{C}$, and $y$ the farthest away 
\begin{tikzpicture}[scale = 0.3, line width = 0.3mm]
\coordinate (aux) at (0,0);
\foreach \i in {0}
    \draw[line cap = round] (aux)--++(1,\i) coordinate (aux);
\end{tikzpicture} from $x$ such that the path from $y$ to $x$ is still a Motzkin path, and
\item $\mathcal{B}$ is the path that remains after removing the path from the first return of the path from the right and the Motzkin path from $y$ to $x$ (including $x$) from the original path.
\end{itemize}

\begin{figure}[h]
\label{canonical}
	\begin{center}
		\begin{tikzpicture}[scale=0.5]
		
		
		\draw (2-2,0) .. controls (3-2,4) .. (4-2,0);
		\draw (4-2,0) .. controls (5-2,4) .. (6-2,0);
		\node at (5,0) {$\cdot$};
		\node at (5.5,0) {$\cdot$};
		\node at (4.5,0) {$\cdot$};
		\draw (6,0) .. controls (7,4) .. (8,1);
		\draw [ultra thick, line cap = round](8,1) to  (9,2);
		\draw [ultra thick, line cap = round](9,2) to  (10,2);

		\node at (5.5+7,2) {$\cdot$};
		\node at (6+7,2) {$\cdot$};
		\node at (6.5+7,2) {$\cdot$};
		\draw (10,2) .. controls (11,5) .. (12,2);	
		\draw (14,2) .. controls (15,5) .. (16,2);	

		\draw [ultra thick, line cap = round](16,2) to  (17,2);
		\draw [ultra thick, line cap = round](17,2) to  (19,0);
		\draw [ultra thick, line cap = round](19,0) to  (20,1);
		\draw (19+1,1) .. controls (20+1,5) .. (21+1,1);
		\draw (21+1,1) .. controls (22+1,5) .. (23+1,1);
		\node at (23.5+1,1) {$\cdot$};
		\node at (24+1,1) {$\cdot$};
		\node at (24.5+1,1) {$\cdot$};
		\draw (25+1,1) .. controls (26+1,5) .. (27+1,1);
			\draw [ultra thick, line cap = round](28,1) to  (29,0);
		
		\node at (9.5,1.5) {$y$};
		\node at (16.5,1.5) {$x$};

		\draw [thick, decorate, decoration={brace, amplitude=10pt, mirror, raise=4pt}] (0cm, -0.5) to
		node[below,yshift=-0.5cm] {$\mathcal{B}_{1}$}
		(8.9cm, -0.5);
		
		\draw [thick, decorate, decoration={brace, amplitude=5pt, mirror, raise=4pt}] (17cm, -0.5) to
		node[below,yshift=-0.5cm] {$\mathcal{B}_{2}$}
		(19cm, -0.5);

		\draw [thick, decorate, decoration={brace, amplitude=10pt, mirror, raise=4pt}] (9.1cm, -0.5) to
		node[below,yshift=-0.5cm] {$\mathcal{A}$}
		(16cm,- 0.5);
		
		\draw [thick, decorate, decoration={brace, amplitude=10pt, mirror, raise=4pt}] (20.1cm, -0.5) to
		node[below,yshift=-0.5cm] {$\mathcal{C}$}
		(28cm, -0.5);
		
		\end{tikzpicture}
		\end	{center}
		
		\caption{Canonical decomposition}

\end{figure}
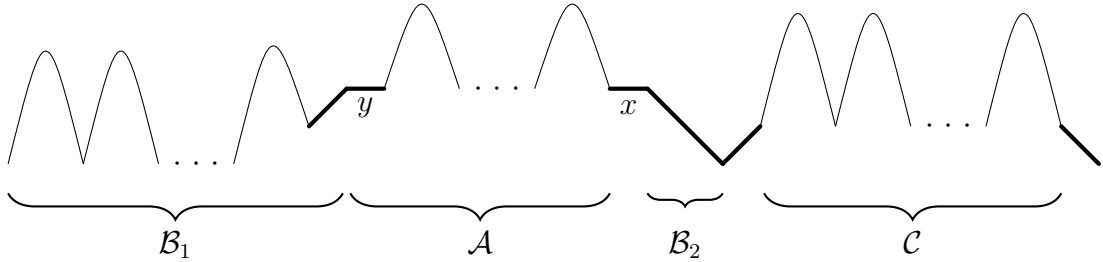

This process is performed recursively and terminates at an empty path. Note that each application of $\Phi$ adds one node and removes one of each type of step. This proves inductively that an S-Motzkin path of length $3n$ maps to a ternary tree with $n$ nodes.

\subsection{Ternary trees to S-Motzkin paths}

The inverse mapping is performed recursively bottom-up as follows. Each node of a ternary tree has three subtrees. Call the paths associated with the left, middle, and right subtrees $\mathcal{A}$, $\mathcal{B}$, and $\mathcal{C}$ respectively. 

Starting at the end nodes, replace each node with
\begin{equation*}
\mathcal{B}_{1}\,\mathcal{A}\,
\begin{tikzpicture}[scale = 0.3, line width = 0.3mm]
\coordinate (aux) at (0,0);
\foreach \i in {0}
    \draw[line cap = round] (aux)--++(1,\i) coordinate (aux);
\end{tikzpicture}\,
\mathcal{B}_{2}\,
\begin{tikzpicture}[scale = 0.3, line width = 0.3mm]
\coordinate (aux) at (0,0);
\foreach \i in {1}
    \draw[line cap = round] (aux)--++(1,\i) coordinate (aux);
\end{tikzpicture}\,
\mathcal{C} \,
\begin{tikzpicture}[scale = 0.3, line width = 0.3mm]
\coordinate (aux) at (0,0);
\foreach \i in {-1}
    \draw[line cap = round] (aux)--++(1,\i) coordinate (aux);
\end{tikzpicture}\,,
\end{equation*}
where $\mathcal{B}_{1}$ is the subpath of $\mathcal{B}$ that starts at $(0, 0)$ and extends to and includes the first occurring 
\begin{tikzpicture}[scale = 0.3, line width = 0.3mm]
\coordinate (aux) at (0,0);
\foreach \i in {1}
    \draw[line cap = round] (aux)--++(1,\i) coordinate (aux);
\end{tikzpicture}
from the right. The path $\mathcal{B}_{2}$ is what remains of $\mathcal{B}$ after removing $\mathcal{B}_{1}$.

This process is continued recursively on each set of end nodes and terminates at the root to produce an \mbox{S-Motzkin} path. Note that for each node three steps are added, and thus a ternary tree with $n$ nodes produces an S-Motzkin path of length $3n$. 

\subsection{Example}

As an example, we map the following S-Motzkin path into a ternary tree. Since the steps are reversible, the inverse mapping can be seen by reading the example in reverse.
Let $\mathcal{M}$ be
\begin{center} 
\begin{tikzpicture}[scale = 0.3, line width = 0.4mm]
\coordinate (aux) at (0,0);
\foreach \i in {0, 1, 0, 1, 0, -1, 1, -1, 0, 1, -1, -1, 0, 1, 0, 1, -1, 0, -1, 1, -1}
    \draw[line cap = round] (aux)--++(1,\i) coordinate (aux);
\addvmargin{1mm}
\end{tikzpicture} .
\end{center}
The canonical decomposition of $\mathcal{M}$ is 
\begin{equation*}
\Phi(\mathcal{M}) = \Big( 
\begin{tikzpicture}[scale = 0.3, line width = 0.4mm]
\coordinate (aux) at (0,0);
\foreach \i in {0, 1, -1}
    \draw[line cap = round] (aux)--++(1,\i) coordinate (aux);
\end{tikzpicture}
\,, \, 
\begin{tikzpicture}[scale = 0.3, line width = 0.4mm]
\coordinate (aux) at (0,0);
\foreach \i in {0, 1, 0, 1, 0, -1, 1, -1, 0, 1, -1, -1, 0, 1, -1}
    \draw[line cap = round] (aux)--++(1,\i) coordinate (aux);
\end{tikzpicture}	
\,, \,
\varnothing
\Big).
\end{equation*}
Hence
\begin{center}
\begin{tikzpicture}
[level distance=12mm,rotate=00,
level 1/.style={sibling distance=40mm},
level 2/.style={sibling distance=20mm},
level 3/.style={sibling distance=20mm}]
\node {$\bullet$}
child {node {
\begin{tikzpicture}[scale = 0.3, line width = 0.4mm]
\coordinate (aux) at (0,0);
\foreach \i in {0, 1, -1}
    \draw[line cap = round] (aux)--++(1,\i) coordinate (aux);
\end{tikzpicture}
}
}
child {node {
\begin{tikzpicture}[scale = 0.3, line width = 0.4mm]
\coordinate (aux) at (0,0);
\foreach \i in {0, 1, 0, 1, 0, -1, 1, -1, 0, 1, -1, -1, 0, 1, -1}
    \draw[line cap = round] (aux)--++(1,\i) coordinate (aux);
\end{tikzpicture}
}
}
child {node {$\varnothing$}};
\end{tikzpicture}
\end{center}
Continuing recursively:
 
$\mathbin{\hphantom{\rightarrow}}$
\begin{minipage}{0.45\textwidth}
\centering
\begin{tikzpicture}
[level distance=10mm,rotate=00,
level 1/.style={sibling distance=10mm},
level 2/.style={sibling distance=15mm},
level 3/.style={sibling distance=15mm}]
\node {$\bullet$}
child {node {$\bullet$}
}
child {node {$\bullet$}
	child{ node{
\begin{tikzpicture}[scale = 0.2, line width = 0.5mm]
\coordinate (aux) at (0,0);
\foreach \i in {0, 1, 0, 1, 0, -1, 1, -1, 0, 1, -1, -1}
    \draw[line cap = round] (aux)--++(1,\i) coordinate (aux);
\addvmargin{1mm}
\end{tikzpicture}
}}
	child{ node{$\varnothing$}}
	child{ node{$\varnothing$}}
}
child[missing] {};
\end{tikzpicture}
\end{minipage}$\rightarrow$
\begin{minipage}{.45\textwidth}
\centering
\begin{tikzpicture}
[level distance=10mm,rotate=00,
level 1/.style={sibling distance=10mm},
level 2/.style={sibling distance=10mm},
level 3/.style={sibling distance=15mm}]
\node {$\bullet$}
child {node {$\bullet$}
}
child {node {$\bullet$}
	child{ node{$\bullet$}
		child{ node{$\varnothing$}}
		child{ node{$\varnothing$}}
		child{ node{
\begin{tikzpicture}[scale = 0.2, line width = 0.5mm]
\coordinate (aux) at (0,0);
\foreach \i in {0, 1, 0, -1, 1, -1, 0, 1, -1}
    \draw[line cap = round] (aux)--++(1,\i) coordinate (aux);
\addvmargin{1mm}
\end{tikzpicture}
}}
}
	child[missing]{}
	child[missing]{}
}
child[missing] {};
\end{tikzpicture}
\end{minipage}

$\rightarrow$
\begin{minipage}{.45\textwidth}
\centering
\begin{tikzpicture}
[level distance=7mm,rotate=00,
level 1/.style={sibling distance=10mm},
level 2/.style={sibling distance=10mm},
level 3/.style={sibling distance=15mm}]
\node {$\bullet$}
child {node {$\bullet$}
}
child {node {$\bullet$}
	child{ node{$\bullet$}
		child[missing]{}
		child[missing]{}
		child{ node{$\bullet$}
			child{ node{
\begin{tikzpicture}[scale = 0.3, line width = 0.5mm]
\coordinate (aux) at (0,0);
\foreach \i in {0, 1, 0, -1, 1, -1}
    \draw[line cap = round] (aux)--++(1,\i) coordinate (aux);
\addvmargin{1mm}
\end{tikzpicture}
}}
			child{ node{$\varnothing$}}	
			child{ node{$\varnothing$}}
}
}
	child[missing]{}
	child[missing]{}
}
child[missing] {};
\end{tikzpicture}
\end{minipage}$\rightarrow$
\begin{minipage}{.45\textwidth}
\centering
\begin{tikzpicture}
[level distance=7mm,rotate=00,
level 1/.style={sibling distance=10mm},
level 2/.style={sibling distance=10mm},
level 3/.style={sibling distance=10mm}]
\node {$\bullet$}
child {node {$\bullet$}
}
child {node {$\bullet$}
	child{ node{$\bullet$}
		child[missing]{}
		child[missing]{}
		child{ node{$\bullet$}
			child{ node{$\bullet$}
				child{ node{$\varnothing$}}
				child{ node{
\begin{tikzpicture}[scale = 0.3, line width = 0.5mm]
\coordinate (aux) at (0,0);
\foreach \i in {0, 1, -1}
    \draw[line cap = round] (aux)--++(1,\i) coordinate (aux);
\addvmargin{1mm}
\end{tikzpicture}
}}
				child{ node{$\varnothing$}}
}
			child[missing]{}	
			child[missing]{}
}
}
	child[missing]{}
	child[missing]{}
}
child[missing] {};
\end{tikzpicture}
\end{minipage}

$\rightarrow$
\begin{minipage}{.45\textwidth}
\centering
\begin{tikzpicture}
[level distance=7mm,rotate=00,
level 1/.style={sibling distance=10mm},
level 2/.style={sibling distance=10mm},
level 3/.style={sibling distance=10mm}]
\node {$\bullet$}
child {node {$\bullet$}
}
child {node {$\bullet$}
	child{ node{$\bullet$}
		child[missing]{}
		child[missing]{}
		child{ node{$\bullet$}
			child{ node{$\bullet$}
				child[missing]{}
				child{ node{$\bullet$}
					child{ node{$\varnothing$}}
					child{ node{$\varnothing$}}
					child{ node{$\varnothing$}}
}
				child[missing]{}
}
			child[missing]{}	
			child[missing]{}
}
}
	child[missing]{}
	child[missing]{}
}
child[missing] {};
\end{tikzpicture}
\end{minipage}$\rightarrow$
\begin{minipage}{.45\textwidth}
\centering
\begin{tikzpicture}
[level distance=7mm,rotate=00,
level 1/.style={sibling distance=10mm},
level 2/.style={sibling distance=10mm},
level 3/.style={sibling distance=10mm}]
\node {$\bullet$}
child {node {$\bullet$}
}
child {node {$\bullet$}
	child{ node{$\bullet$}
		child[missing]{}
		child[missing]{}
		child{ node{$\bullet$}
			child{ node{$\bullet$}
				child[missing]{}
				child{ node{$\bullet$}}
				child[missing]{}
}
			child[missing]{}	
			child[missing]{}
}
}
	child[missing]{}
	child[missing]{}
}
child[missing] {};
\end{tikzpicture}
\end{minipage}

\begin{table}[h!]
\caption{Bijection for $n=3$}
\setlength{\tabcolsep}{5mm} 
\def\arraystretch{1.25} 
\centering
\begin{tabular}{|c|c|}
  \hline
  S-Motzkin path   &   Ternary tree 
  \\ \hline
  
\begin{tikzpicture}[scale = 0.4, line width = 0.5mm]
\coordinate (aux) at (0,0);
\foreach \i in {0, 1, -1, 0, 1, -1, 0, 1, -1}
    \draw[line cap = round] (aux)--++(1,\i) coordinate (aux);
\addvmargin{1mm}
\end{tikzpicture}  &    
\begin{tikzpicture}
[level distance=5mm,rotate=00,
level 1/.style={sibling distance=7mm},
level 2/.style={sibling distance=7mm},
level 3/.style={sibling distance=7mm}]
\node {$\bullet$}
child {node {$\bullet$} 
	child{ node {$\bullet$}
			}
	child[missing]{ 
			}
	child[missing]{
			}
}
child[missing] {
}
child[missing] {
};
\end{tikzpicture}
  
  \\ \hline

  \begin{tikzpicture}[scale = 0.4, line width = 0.5mm]
\coordinate (aux) at (0,0);
\foreach \i in {0, 1, -1, 0, 1, 0, -1, 1, -1}
    \draw[line cap = round] (aux)--++(1,\i) coordinate (aux);
\addvmargin{1mm}
\end{tikzpicture}
   &   \begin{tikzpicture}
[level distance=7mm,rotate=00,
level 1/.style={sibling distance=7mm},
level 2/.style={sibling distance=7mm},
level 3/.style={sibling distance=7mm}]
\node {$\bullet$}
child[missing] {}
child {node {$\bullet$}
	child {node {$\bullet$}
			}
	child[missing] {}
	child[missing] {}
}
child[missing] {};
\end{tikzpicture}
   
  \\ \hline

\begin{tikzpicture}[scale = 0.4, line width = 0.5mm]
\coordinate (aux) at (0,0);
\foreach \i in {0, 1, -1, 0, 1, 0, 1, -1, -1}
    \draw[line cap = round] (aux)--++(1,\i) coordinate (aux);
\addvmargin{1mm}
\end{tikzpicture}   &  
\begin{tikzpicture}
[level distance=5mm,rotate=00,
level 1/.style={sibling distance=7mm},
level 2/.style={sibling distance=7mm},
level 3/.style={sibling distance=7mm}]
\node {$\bullet$}
child {node {$\bullet$} }
child[missing] {}
child {node {$\bullet$}};
\end{tikzpicture}

  \\ \hline

\begin{tikzpicture}[scale = 0.4, line width = 0.5mm]
\coordinate (aux) at (0,0);
\foreach \i in {0, 1, 0, -1, 1, -1, 0, 1, -1}
    \draw[line cap = round] (aux)--++(1,\i) coordinate (aux);
\addvmargin{1mm}
\end{tikzpicture}
   &  
\begin{tikzpicture}
[level distance=7mm,rotate=00,
level 1/.style={sibling distance=7mm},
level 2/.style={sibling distance=7mm},
level 3/.style={sibling distance=7mm}]
\node {$\bullet$}
child {node {$\bullet$} 
	child[missing]{ }
	child{ node {$\bullet$} 
			}
	child[missing]{}
}
child[missing] {}
child[missing] {};
\end{tikzpicture}

  \\ \hline

\begin{tikzpicture}[scale = 0.4, line width = 0.5mm]
\coordinate (aux) at (0,0);
\foreach \i in {0, 1, 0, -1, 1, 0, -1, 1, -1}
    \draw[line cap = round] (aux)--++(1,\i) coordinate (aux);
\addvmargin{1mm}
\end{tikzpicture}
   &  
\begin{tikzpicture}
[level distance=7mm,rotate=00,
level 1/.style={sibling distance=7mm},
level 2/.style={sibling distance=7mm},
level 3/.style={sibling distance=7mm}]
\node {$\bullet$}
child[missing] {}
child {node {$\bullet$}
	child[missing] {}
	child {node {$\bullet$}
			}
	child[missing] {}
}
child[missing] {};
\end{tikzpicture}

  \\ \hline

\begin{tikzpicture}[scale = 0.4, line width = 0.5mm]
\coordinate (aux) at (0,0);
\foreach \i in {0, 1, 0, -1, 1, 0, 1, -1, -1}
    \draw[line cap = round] (aux)--++(1,\i) coordinate (aux);
\addvmargin{1mm}
\end{tikzpicture}
   &  
\begin{tikzpicture}
[level distance=7mm,rotate=00,
level 1/.style={sibling distance=7mm},
level 2/.style={sibling distance=7mm},
level 3/.style={sibling distance=7mm}]
\node {$\bullet$}
child[missing] {}
child {node {$\bullet$}}
child {node {$\bullet$}};
\end{tikzpicture}

  \\ \hline

\begin{tikzpicture}[scale = 0.4, line width = 0.5mm]
\coordinate (aux) at (0,0);
\foreach \i in {0, 1, 0, 1, -1, 0, 1, -1, -1}
    \draw[line cap = round] (aux)--++(1,\i) coordinate (aux);
\addvmargin{1mm}
\end{tikzpicture}  &  
\begin{tikzpicture}
[level distance=5mm,rotate=00,
level 1/.style={sibling distance=7mm},
level 2/.style={sibling distance=7mm},
level 3/.style={sibling distance=7mm}]
\node {$\bullet$}
child[missing] {}
child[missing] {}
child {node {$\bullet$}
	child {node {$\bullet$}
			}
	child[missing] {}
	child[missing] {}
};
\end{tikzpicture}

  \\ \hline

\begin{tikzpicture}[scale = 0.4, line width = 0.5mm]
\coordinate (aux) at (0,0);
\foreach \i in {0, 1, 0, 1, -1, 0, -1, 1, -1}
    \draw[line cap = round] (aux)--++(1,\i) coordinate (aux);
\addvmargin{1mm}
\end{tikzpicture}   &  
\begin{tikzpicture}
[level distance=7mm,rotate=00,
level 1/.style={sibling distance=7mm},
level 2/.style={sibling distance=7mm},
level 3/.style={sibling distance=7mm}]
\node {$\bullet$}
child {node {$\bullet$}}
child {node {$\bullet$}}
child[missing] {};
\end{tikzpicture}

  \\ \hline

\begin{tikzpicture}[scale = 0.4, line width = 0.5mm]
\coordinate (aux) at (0,0);
\foreach \i in {0, 1, 0, 1, -1, -1, 0, 1, -1}
    \draw[line cap = round] (aux)--++(1,\i) coordinate (aux);
\addvmargin{1mm}
\end{tikzpicture}   &  
\begin{tikzpicture}
[level distance=5mm,rotate=00,
level 1/.style={sibling distance=7mm},
level 2/.style={sibling distance=7mm},
level 3/.style={sibling distance=7mm}]
\node {$\bullet$}
child {node {$\bullet$} 
	child[missing]{ }
	child[missing]{ }
	child{ node {$\bullet$}
			}
}
child[missing] {}
child[missing] {};
\end{tikzpicture}

  \\ \hline

\begin{tikzpicture}[scale = 0.4, line width = 0.5mm]
\coordinate (aux) at (0,0);
\foreach \i in {0, 1, 0, 1, 0, -1, -1, 1, -1}
    \draw[line cap = round] (aux)--++(1,\i) coordinate (aux);
\addvmargin{1mm}
\end{tikzpicture}   &  
\begin{tikzpicture}
[level distance=7mm,rotate=00,
level 1/.style={sibling distance=7mm},
level 2/.style={sibling distance=7mm},
level 3/.style={sibling distance=7mm}]
\node {$\bullet$}
child[missing] {}
child {node {$\bullet$}
	child[missing] {}
	child[missing] {}
	child {node {$\bullet$}
			}
}
child[missing] {};
\end{tikzpicture}

  \\ \hline

\begin{tikzpicture}[scale = 0.4, line width = 0.5mm]
\coordinate (aux) at (0,0);
\foreach \i in {0, 1, 0, 1, 0, -1, 1, -1, -1}
    \draw[line cap = round] (aux)--++(1,\i) coordinate (aux);
\addvmargin{1mm}
\end{tikzpicture}  &  
\begin{tikzpicture}
[level distance=7mm,rotate=00,
level 1/.style={sibling distance=7mm},
level 2/.style={sibling distance=7mm},
level 3/.style={sibling distance=7mm}]
\node {$\bullet$}
child[missing] {}
child[missing] {}
child {node {$\bullet$}
	child[missing] {}
	child {node {$\bullet$}
			}
	child[missing] {}
};
\end{tikzpicture}

  \\ \hline

\begin{tikzpicture}[scale = 0.4, line width = 0.5mm]
\coordinate (aux) at (0,0);
\foreach \i in {0, 1, 0, 1, 0, 1, -1, -1, -1}
    \draw[line cap = round] (aux)--++(1,\i) coordinate (aux);
\addvmargin{1mm}
\end{tikzpicture}  &  
\begin{tikzpicture}
[level distance=5mm,rotate=00,
level 1/.style={sibling distance=7mm},
level 2/.style={sibling distance=7mm},
level 3/.style={sibling distance=7mm}]
\node {$\bullet$}
child[missing] {}
child[missing] {}
child {node {$\bullet$}
	child[missing] {}
	child[missing] {}
	child {node {$\bullet$}
			}
};
\end{tikzpicture}

  \\ \hline

  \end{tabular}
\end{table}

\clearpage

\section{Generalization}
This bijection can be generalized to be between $t$-ary trees and the subclass of Motzkin paths with $(t-2)n$ 
\,\begin{tikzpicture}[scale = 0.3, line width = 0.3mm]
\coordinate (aux) at (0,0);
\foreach \i in {0}
    \draw[line cap = round] (aux)--++(1,\i) coordinate (aux);
\end{tikzpicture}\, steps and $n$ of each of the other steps such that 
\begin{itemize}
\item The initial $t-2$ steps must be of the form 
\begin{tikzpicture}[scale = 0.3, line width = 0.3mm]
\coordinate (aux) at (0,0);
\foreach \i in {0}
    \draw[line cap = round] (aux)--++(1,\i) coordinate (aux);
\end{tikzpicture}\,, 
\item between every two 
\begin{tikzpicture}[scale = 0.3, line width = 0.3mm]
\coordinate (aux) at (0,0);
\foreach \i in {1}
    \draw[line cap = round] (aux)--++(1,\i) coordinate (aux);
\end{tikzpicture} 
there are exactly $t-2$ steps of the form 
\begin{tikzpicture}[scale = 0.3, line width = 0.3mm]
\coordinate (aux) at (0,0);
\foreach \i in {0}
    \draw[line cap = round] (aux)--++(1,\i) coordinate (aux);
\end{tikzpicture}\,, 
\item the $k$-th occurring 
\begin{tikzpicture}[scale = 0.3, line width = 0.3mm]
\coordinate (aux) at (0,0);
\foreach \i in {-1}
    \draw[line cap = round] (aux)--++(1,\i) coordinate (aux);
\end{tikzpicture}
must occur after at least $k$ occurrences of $t-2$ steps of the form
\begin{tikzpicture}[scale = 0.3, line width = 0.3mm]
\coordinate (aux) at (0,0);
\foreach \i in {0}
    \draw[line cap = round] (aux)--++(1,\i) coordinate (aux);
\end{tikzpicture}
and one step of the form
\begin{tikzpicture}[scale = 0.3, line width = 0.3mm]
\coordinate (aux) at (0,0);
\foreach \i in {1}
    \draw[line cap = round] (aux)--++(1,\i) coordinate (aux);
\end{tikzpicture}\,.
\end{itemize}
\phantom{.}

\textbf{Acknowledgement.} The authors thank Stephan Wagner for valuable feedback. 

\bibliography{bijection_2}

\end{document}